\documentclass[12pt,leqno]{amsart}
\usepackage{amssymb,amsthm,verbatim,enumerate,hyperref,ifthen}
\usepackage[mathscr]{eucal}
\usepackage[utf8]{inputenc}
\usepackage[T1]{fontenc}
\usepackage{graphics}
\usepackage{graphicx}
\def\N{\mathbb{N}}
\def\R{\mathbb{R}}

\def\Z{\mathbb{Z}}
\def\A{\mathscr{A}}

\def\dist{\mathop{\mbox{\rm dist}}}

\def\epsilon{\varepsilon}
\def\O{\Omega}
\def\o{\omega}
\def\phi{\varphi}
\def\eps{\varepsilon}

\newtheorem{theorem}{Theorem}

\newtheorem*{theorem*}{Theorem}
\def\Thm#1#2{\ifthenelse{\equal{#1}{*}}{\begin{theorem*}#2\end{theorem*}}
  {\begin{theorem}\label{T#1}#2\end{theorem}}}
\newtheorem{Atheorem}{Theorem}

\def\THM#1#2{\begin{Atheorem}\label{T#1}#2\end{Atheorem}}
\def\thm#1{Theorem~\ref{T#1}}
\newtheorem{proposition}[theorem]{Proposition}
\newtheorem*{proposition*}{Proposition}
\def\Prp#1#2{\ifthenelse{\equal{#1}{*}}{\begin{proposition*}#2\end{proposition*}}
             {\begin{proposition}\label{P#1}#2\end{proposition}}}
\def\prp#1{Proposition~\ref{P#1}}
\newtheorem{corollary}[theorem]{Corollary}
\newtheorem*{corollary*}{Corollary}
\def\Cor#1#2{\ifthenelse{\equal{#1}{*}}{\begin{corollary*}#2\end{corollary*}}
             {\begin{corollary}\label{C#1}#2\end{corollary}}}
\def\cor#1{Corollary~\ref{C#1}}
\newtheorem{lemma}[theorem]{Lemma}
\newtheorem*{lemma*}{Lemma}
\def\Lem#1#2{\ifthenelse{\equal{#1}{*}}{\begin{lemma*}#2\end{lemma*}}
             {\begin{lemma}\label{L#1}#2\end{lemma}}}

\newtheorem{example}[theorem]{Example}
\newtheorem*{example*}{Example}
\long\def\Exa#1#2{\ifthenelse{\equal{#1}{*}}{\begin{example*}\rm #2\end{example*}}
            {\begin{example}\label{Ex#1}\rm #2\end{example}}}

\newtheorem{problem}[subsection]{Problem}

\theoremstyle{definition}
\newtheorem{definition}[theorem]{Definition}
\newtheorem*{definition*}{Definition}
\def\Defi#1#2{\ifthenelse{\equal{#1}{*}}{\begin{definition*}#2\end{definition*}}
      {\begin{definition}\label{D#1}#2\end{definition}}}

\newtheorem{remark}[theorem]{Remark}
\newtheorem*{remark*}{Remark}
\def\Rem#1#2{\ifthenelse{\equal{#1}{*}}{\begin{remark*}#2\end{remark*}}
             {\begin{remark}\label{R#1}#2\end{remark}}}
\def\rem#1{Remark~\ref{R#1}}
\def\eq#1{{\rm(\ref{E#1})}}

\def\Eq#1#2{\ifthenelse{\equal{#1}{*}}
  {\begin{equation*}\begin{aligned}{}#2\end{aligned}\end{equation*}}
  {\begin{equation}\begin{aligned}{}\label{E#1}#2\end{aligned}\end{equation}}}

\begin{document}
\begin{flushright}
Math. Inequal. Appl. \textbf{16}(2) (2013), 507--526. \\
\href{http://dx.doi.org/10.7153/mia-16-37}{doi: 10.7153/mia-16-37} \\[1cm]
\end{flushright}

\title[Approximate Hermite--Hadamard type inequalities]
{Approximate Hermite--Hadamard type inequalities for approximately convex functions}
\author[J.\ Mak\'o]{Judit Mak\'o}
\author[Zs. P\'ales]{Zsolt P\'ales}
\address{Institute of Mathematics, University of Debrecen,
H-4010 Debrecen, Pf.\ 12, Hungary}
\email{\{makoj,pales\}@science.unideb.hu} 
\thanks{This research has been supported by the Hungarian
Scientific Research Fund (OTKA) Grant NK81402 and by the
TÁMOP-4.2.1/B-09/1/KONV-2010-0007, TÁMOP-4.2.2/B-10/1-2010-0024 projects. 
These projects are co-financed by the European Union and the European Social Fund.}

\begin{abstract}
In this paper, approximate lower and upper Hermite--Hadamard type inequalities are obtained for
functions that are approximately convex with respect to a given Chebyshev system.
\end{abstract}

\maketitle

\section{Introduction}

Throughout this paper $\R$, $\R_+$, $\N$ and $\Z$ denote the sets of
real, nonnegative real, natural and integer numbers respectively.
Given a nonempty open real interval $I$, denote by $\Delta (I)$ and $\Delta^\circ(I)$ the sets
\Eq{*}{
\{(x,y)\in I\times I\mid x\leq y\} \qquad\mbox{and}\qquad\{(x,y)\in I\times I\mid x<y\},
}
respectively. Given a nonempty open real interval $I$, denote by $\Delta (I)$ and $\Delta^\circ(I)$ the sets
\Eq{*}{
\{(x,y)\in I\times I\mid x\leq y\} \qquad\mbox{and}\qquad\{(x,y)\in I\times I\mid x<y\},
}
respectively.
We say that a pair $(\o_0,\o_1)$ is a \textit{Chebyshev system} over $I$, if
$\o_0,\o_1:I\to\R$ are continuous functions and
\Eq{O}{
\O(x,y):=\left| \begin{array}{ccc}
\o_0(x) & \o_0(y)  \\
\o_1(x) & \o_1(y)  \\
\end{array} \right|>0 \qquad((x,y)\in\Delta^\circ(I)).
}
One can easily see, that if $\o_0$ is a positive function, then \eq{O} holds if and only if $\o_1/\o_0$ is
strictly increasing on $I$. In this latter case, $(\o_0,\o_1)$ will be called a \textit{positive Chebyshev system}
over $I$. On the other hand, we can always assume that $\o_0$ is a positive function, because
for every Chebyshev system $(\o_0,\o_1)$, there exists $\alpha,\beta\in\R$ such that $\alpha\o_0+\beta\o_1>0$
(cf.\ \cite{Bes08b}, \cite{BesPal03}). In the sequel, for fixed $x,y\in I$, the partial functions
$u\mapsto\O(u,y)$ and $u\mapsto\O(x,u)$ will be denoted by $\O(\cdot,y)$ and $\O(x,\cdot)$, respectively.
An important property of Chebyshev systems is that for every two pairs $(x,\xi),(y,\eta)\in I\times\R$ with $x\neq y$ the
function $\o$ defined as
\Eq{*}{
  \o:=\xi\frac{\O(\cdot,y)}{\O(x,y)}+\eta\frac{\O(x,\cdot)}{\O(x,y)}
}
is the unique linear combination of $\o_0$ and $\o_1$ such that $\o(x)=\xi$ and $\o(y)=\eta$ hold.

Given a positive Chebyshev system $(\o_0,\o_1)$ over $I$ and a proper subinterval $J$ of $I$,
a function $f:J\to\R$ is called \textit{$(\o_0,\o_1)$-convex on $J$} if, for all $x<u<y$ from $J$,
\Eq{OC1}{
\left| \begin{array}{ccc}
f(x)      &f(u)         &f(y) \\
\o_0(x) & \o_0(u)   &\o_0(y)\\
\o_1(x) & \o_1(u)   &\o_1(y)
\end{array} \right|\geq0,
}
or equivalently,
\Eq{OC2}{
  f(u)\leq \frac{\O(u,y)}{\O(x,y)}f(x)+\frac{\O(x,u)}{\O(x,y)}f(y).
}
If \eq{OC1} holds with strict inequality sign ``$>$'', then $f$ is said to be
\textit{strictly $(\o_0,\o_1)$-convex on $J$}.

The integral average of any standard convex function $f : I \to \R$ can be estimated from the midpoint and
the endpoints of the domain as follows:
\Eq{nrsl}{
     f\Big(\frac{x+y}{2}\Big)\leq\int\limits_0^1 f\big(tx+(1-t)y\big)dt\leq\frac{f(x)+f(y)}{2}\qquad(x,y\in I).
}
This is the well known Hermite--Hadamard type inequality. The above implication was discovered by
Hadamard \cite{Had1893}. (See also \cite{MitLac85}, \cite{Kuc85}, and \cite{NicPer06},
\cite{DraPea00}, \cite{NicPer03}, \cite{NicPer06}, \cite{NikRieSah07} for a historical account.)
In \cite{BesPal06} and \cite{Bes08b}, the authors established the following connections between
$(\o_0,\o_1)$-convexity and Hermite--Hadamard type inequality.

\THM{E}{Let $(\o_0,\o_1)$ be a positive Chebyshev system on an open interval $I$ and
let $\rho: I \to\R$ be a positive integrable function. Define, for all elements $x<y$ of $I$,
the functions $\xi(x,y)$, $c(x,y)$, $c_1(x,y)$ and $c_2(x,y)$ by the formulas
\Eq{0kc}{
&\xi(x,y)=\Big(\frac{\o_1}{\o_0}\Big)^{-1}\Big(\frac{\int_x^y\o_1\rho}{\int_x^y\o_0\rho}\Big)
\quad\mbox{and}\quad c(x,y)=\frac{\int_x^y\o_0\rho}{\o_0(\xi(x,y))},\\
&c_1(x,y)=\frac{1}{\O(x,y)}\left| \begin{array}{ccc}
\int_x^y\o_0\rho &\o_0(y) \\[2mm]
\int_x^y\o_1\rho & \o_1(y)
\end{array} \right|\quad \mbox{and}\quad
c_2(x,y)=\frac{1}{\O(x,y)}\left| \begin{array}{ccc}
\o_0(x) &\int_x^y\o_0\rho \\[2mm]
\o_0(y) &\int_x^y\o_1\rho
\end{array} \right|.
}
If a function $f : I \to \R$ is $(\o_0,\o_1)$-convex, then for all elements
$x<y$ of $I$, it satisfies the inequality
\Eq{EH}{
c(x,y)f(\xi(x,y))\leq \int\limits_x^yf\rho\leq c_1(x,y)f(x)+c_2(x,y)f(y).
}
}

In \thm{3} and \thm{5} below, these results will be generalized to the context 
of approximate $(\o_0,\o_1)$-convex, i.e., to the case when $f$ satisfies an inequality analogous to
\eq{OC2} whose right hand side involves also an error term.

Let $X$ be a real linear space and $D\subset X$ be a convex subset.
In order to describe the old and new results about the connection of
approximate Jensen convexity and the approximate Hermite--Hadamard inequality
with variable error terms, we need to introduce the following terminology.

For a function $f:D\to\R$, we say that $f$ is \textit{hemi-}$P$, if, for all $x,y\in D$, the mapping
\Eq{0f}
{
t\mapsto f((1-t)x+ty) \qquad (t\in [0,1])
}
has property $P$. For example $f$ is hemiintegrable, if for all $x,y\in D$ the mapping defined 
by \eq{0f} is integrable. Analogously, we say that a function $h:(D-D)\to\R$ is
\textit{radially}-$P$, if for all $u\in D-D$, the mapping
\Eq{*}
{
t\mapsto h(tu) \qquad (t\in [0,1])
}
has property $P$ on $[0,1]$.

In \cite{HazPal09}, H\'azy and second author of this paper established a connection between an
approximate lower
Hermite--Hadamard type inequality and an approximate Jensen type inequality by proving the following
result.

\THM{C}{Let $\alpha:(D-D)\to\R_+$ be a nonnegative radially Lebesgue integrable even function.
Assume that $f:D\to\R$ is hemi-Lebesgue integrable and approximately Jensen convex in the sense of
\Eq{0J}{
f\Big(\frac{x+y}{2}\Big)\leq\frac{f(x)+f(y)}{2}+\alpha(x-y) \qquad (x,y\in D).
}
Then $f$ also satisfies the approximate lower Hermite--Hadamard inequality
\Eq{0lH}{
f\left(\frac{x+y}{2}\right)\leq \int\limits_0^1 f\big(tx+(1-t)y\big)dt +
  \int\limits_0^1 \alpha(t(x-y))dt \qquad(x,y\in D).
}}

In \cite{MakPal11d} (cf.\ \cite{TabTab09b}, \cite{TabTabZol10a}) the authors established the
connections 
between an approximate upper Hermite--Hadamard
type inequality and an approximate Jensen type inequality as stated in the following theorem.

\THM{A}{Let $\alpha:(D-D)\to\R_+$ be a nonnegative radially Lebesgue integrable even function 
and $\rho:[0,1]\to\R_+$ be a nonnegative Lebesgue integrable function with $\int_0^1\rho=1$.
Assume that $f:D\to\R$ is hemiintegrable on $D$ and satisfies the approximate Jensen inequality
\eq{0J}.
Then, for $x,y\in D$, $f$ also satisfies the approximate upper Hermite--Hadamard inequality
\Eq{0H}{
\int_0^1f\big(tx+(1-t)y\big)\rho(t)dt
\leq\lambda f(x)+(1-\lambda)f(y) +\sum_{n=0}^\infty \frac{1}{2^n}
   \int_{0}^1\alpha\big(2d_{\Z}(2^{n}t)(x-y)\big)\rho(t)dt,
}
where $\lambda:=\int_0^1 t\rho(t)dt$ and, for $s\in\R$, $d_\Z(s):=\dist(s,\Z)=\inf\{|s-k|:k\in\Z\}$.
}

In \thm{4} and \thm{6} below these results will be generalized and extended to the setting of
$(\o_0,\o_1)$-convexity.

\section{From approximate $(\o_0,\o_1)$-convexity to approximate lower Hermite--Hadamard inequality}
\setcounter{equation}{0}

In this section we will investigate the implication between an $(\o_0,\o_1)$-convexity
type inequality and a lower Hermite--Hadamard inequality. Consider the following basic assumptions.

\begin{enumerate}[({A}1)]
\item $(T,\A,\mu)$ is a measure space.
\item $\Lambda:T\times \Delta^\circ(I)\to\R_+$ is $\mu$-integrable in its first variable.
\item $M:T\times \Delta^\circ(I)\to\R$ is $\A$-measurable in its first variable
and for all $t\in T$, the map $(x,y)\mapsto M(t,x,y)$ is a two-variable mean on $I$.
$M_0:\Delta^\circ(I)\to I$ is a strict mean such that
\Eq{2fi+}{
   \mu\{t\in T\mid \Lambda(t,x,y)>0,\,M(t,x,y)\neq M_0(x,y)\}>0
        \qquad\mbox{if}\quad (x,y)\in \Delta^{\circ}(I).}
\item There exist an $(\o_0,\o_1)$-Chebyshev system on $I$ such that $\o_0$ is positive.
Furthermore, for $i\in\{0,1\}$,
\Eq{2=}{
\o_i(M_0(x,y))=\int\limits_T\Lambda(t,x,y)\o_i(M(t,x,y))d\mu(t) \qquad ((x,y)\in\Delta^\circ(I)).
}
\end{enumerate}

For all $(x,y)\in \Delta^{\circ}(I)$, denote
\Eq{*}{
&T'_{x,y}:=\{t\in T\mid \Lambda(t,x,y)>0,\, M(t,x,y)<M_0(x,y)\},\\
&T''_{x,y}:=\{t\in T\mid\Lambda(t,x,y)>0,\, M(t,x,y)\geq M_0(x,y)\}.
}
Observe that, for all $(x,y)\in\Delta^{\circ}(I)$, $T'_{x,y}$ and $T''_{x,y}$ are in $\A$,
moreover, by \eq{2fi+}, the $\mu$-measure of $T'_{x,y}\cup T''_{x,y}$ is positive.
Define, for all $(x,y)\in\Delta^{\circ}(I)$, $i\in\{0,1\}$,
\Eq{Si}{
S_i'(x,y)=\int\limits_{T'_{x,y}}\Lambda(t,x,y)\o_i&(M(t,x,y))d\mu(t)\quad
\mbox{and}\quad S_i''(x,y)=\int\limits_{T''_{x,y}}\Lambda(t,x,y)\o_i(M(t,x,y))d\mu(t).
}
The following proposition describes the properties of these sets and numbers.

\Prp{TS}{
If (A1)--(A4) hold, then for all $(x,y)\in \Delta^{\circ}(I)$,
\Eq{TS}{
S_i'(x,y)+S_i''(x,y)=\o_i(M_0(x,y)) \qquad (i\in\{0,1\}).
}
Furthermore, the $\mu$-measure of the sets $T'_{x,y}$ and $T''_{x,y}$
is positive.}

\begin{proof}
Let $(x,y)\in\Delta^{\circ}(I)$. \eq{2=} implies that
\Eq{*}{
\o_i(M_0(x,y))&=\int\limits_{T}\Lambda(t,x,y)\o_i(M(t,x,y))d\mu(t)\\
              &=\int\limits_{\{t\in T\mid \Lambda(t,x,y)>0\}}\Lambda(t,x,y)\o_i(M(t,x,y))d\mu(t)
               =S_i'(x,y)+S_i''(x,y),
}
for $i\in\{0,1\}$. Hence \eq{TS} holds.
To prove the positivity of $\mu(T'_{x,y})$ and $\mu(T''_{x,y})$, assume that $\mu(T'_{x,y})=0$.
Then $S_i'(x,y)=0$ and, in view of \eq{2fi+}, it follows that $\mu(T''_{x,y})>0$.
Thus, by \eq{TS}, we have that
\Eq{*}{
\o_i(M_0(x,y))&=S_i'(x,y)+S_i''(x,y)=S_i''(x,y)=\int\limits_{T''_{x,y}}\Lambda(t,x,y)\o_i(M(t,x,
y))d\mu(t)
}
for $i\in\{0,1\}$. Dividing the above identities by each other and using also the positivity of
$\o_0$, we get that
\Eq{*}{
\frac{\int\limits_{T''_{x,y}}\Lambda(t,x,y)\o_1(M(t,x,y))d\mu(t)}
{\int\limits_{T''_{x,y}}\Lambda(t,x,y)\o_0(M(t,x,y))d\mu(t)}
   =\frac{\o_1(M_0(x,y))}{\o_0(M_0(x,y))}.
}
Rearranging this equality, we obtain that
\Eq{*}{
\int\limits_{T''_{x,y}}\Lambda(t,x,y)\O(M_0(x,y),M(t,x,y))d\mu(t)=0.
}
Hence,
\Eq{*}{
\int\limits_{\{t\in T\mid\Lambda(t,x,y)>0,\, M(t,x,y)>
M_0(x,y)\}}\Lambda(t,x,y)\O(M_0(x,y),M(t,x,y))d\mu(t)=0.
}
On the other hand, for all $t\in T$ with $M(t,x,y)> M_0(x,y)$, we have that
$\O(M_0(x,y),M(t,x,y))>0$ and,
by \eq{2fi+}, $\mu(\{t\in T\mid\Lambda(t,x,y)>0,\, M(t,x,y)> M_0(x,y)\})>0$.
This yields that
\Eq{*}{
\int\limits_{\{t\in T\mid\Lambda(t,x,y)>0,\, M(t,x,y)>
M_0(x,y)\}}\Lambda(t,x,y)\O(M_0(x,y),M(t,x,y))d\mu(t)>0,}
which is a contradiction.

The proof for the case when $\mu(T''_{x,y})=0$ is analogous.
\end{proof}

One of the main result of this paper is established in the following theorem.

\Thm{3}{Assume that (A1)--(A4) hold.
Let $f:I\to\R$ be a locally upper bounded Borel measurable solution of the approximate
$(\o_0,\o_1)$-convexity
type functional inequality
\Eq{J}{
f(u)\leq \frac{\O(u,y)}{\O(x,y)}f(x)+\frac{\O(x,u)}{\O(x,y)}f(y)+\eps_{x,y}(u)\qquad (u\in [x,y]),}
where for all $(x,y)\in \Delta^{\circ}(I)$ and $u\in]x,y[$, the function
$(v,w)\mapsto\eps_{v,w}(u)$ is bounded and Borel measurable for $(v,w)\in[x,u]\times[u,y]$. Then $f$
also
satisfies the approximate lower Hermite--Hadamard type inequality
\Eq{lH}{
  f(M_0(x,y))\leq \int\limits_T \Lambda(t,x,y)f\big(M(t,x,y)\big)d\mu(t) +
  E(x,y) \qquad((x,y)\in \Delta(I)),
}
where $E:\Delta^\circ(I)\to\R$ is defined by the following way
\Eq{2E}{
&E(x,y)\\&=\frac{\int\limits_{T'_{x,y}}\int\limits_{T''_{x,y}}
          \Lambda(t',x,y)\Lambda(t'',x,y)\O(M(t',x,y),M(t'',x,y))
          \eps_{M(t',x,y),M(t'',x,y)}(M_0(x,y))d\mu(t'')d\mu(t')}
          {\int\limits_{T'_{x,y}}\int\limits_{T''_{x,y}}
          \Lambda(t',x,y)\Lambda(t'',x,y)\O(M(t',x,y),M(t'',x,y))d\mu(t'')d\mu(t')}.
}}

\Rem{T3}{In the above theorem, the regularity condition for $f$ can be relaxed if the error function
$\eps_{x,y}$ enjoys boundedness and continuity properties. For instance, if $\eps_{x,y}$ is bounded
on $[x,y]$
for some $(x,y)\in\Delta^\circ(I)$, then \eq{J} implies that $f$ is bounded on $[x,y]$. Similarly,
if
$\limsup_{u\to x+0}\eps_{x,y}(u)=0$ for some $(x,y)\in\Delta^\circ(I)$, then \eq{J} implies that $f$
is
upper semicontinuous at $x$ from the right.}

\begin{proof} Let $(x,y)\in\Delta^{\circ}(I)$. Substituting in \eq{J} $x$ by $M(t',x,y)$ and $y$ by
$M(t'',x,y)$,
and $u$ by $M_0(x,y)$, where $t'\in T'_{x,y}$ and $t''\in T''_{x,y}$, we get that
\Eq{*}{
\O(M(t',x,y),M(t'',x,y))f(M_0(x,y))
          &\leq \O(M_0(x,y),M(t'',x,y))f(M(t',x,y))\\
          &+\O(M(t',x,y),M_0(x,y))f(M(t'',x,y))\\
          &+\O(M(t',x,y),M(t'',x,y))\eps_{M(t',x,y),M(t'',x,y)}(M_0(x,y)).
}
Multiplying this inequality by $\Lambda(t',x,y)\Lambda(t'',x,y)$ and integrating with respect to
$\mu\times\mu$ on $T'_{x,y}\times T''_{x,y}$, we get that
\Eq{F0}{
\int\limits_{T'_{x,y}}&\int\limits_{T''_{x,y}}
\Lambda(t',x,y)\Lambda(t'',x,y)\O(M(t',x,y),M(t'',x,y))d\mu(t'')d\mu(t')f(M_0(x,y))\\
&\leq \int\limits_{T'_{x,y}}\int\limits_{T''_{x,y}}
\Lambda(t',x,y)\Lambda(t'',x,y)\O(M_0(x,y),M(t'',x,y))f(M(t',x,y))d\mu(t'')d\mu(t')\\
&+\int\limits_{T'_{x,y}}\int\limits_{T''_{x,y}}
\Lambda(t',x,y)\Lambda(t'',x,y)\O(M(t',x,y),M_0(x,y))f(M(t'',x,y))d\mu(t'')d\mu(t')\\
&+\int\limits_{T'_{x,y}}\int\limits_{T''_{x,y}}\!\!\!\!
\Lambda(t',x,y)\Lambda(t'',x,y)\O(M(t',x,y),M(t'',x,y))
    \eps_{M(t',x,y),M(t'',x,y)}(M_0(x,y))d\mu(t'')d\mu(t').
}
Applying Fubini's theorem and the notation of \eq{Si}, we get that
\Eq{F1}{
\int\limits_{T'_{x,y}}\int\limits_{T''_{x,y}}&\Lambda(t',x,y)\Lambda(t'',x,y)\O(M(t',x,y),M(t'',x,
y))d\mu(t'')d\mu(t')\\
&=\big(S_0'(x,y)S_1''(x,y)-S_1'(x,y)S_0''(x,y)\big).
}
Observe that $(S_0'(x,y)S_1''(x,y)-S_1'(x,y)S_0''(x,y))$ is positive. Indeed,
by the definition of the Chebyshev-system, we have, for all $(t',t'')\in T'_{x,y}\times T''_{x,y}$,
\Eq{*}{
\O(M(t',x,y),M(t'',x,y))>0.
}
By \prp{TS}, the measure of $T'_{x,y}\times T''_{x,y}$ is positive. Hence the left hand side of
\eq{F1} is positive. Using the identity \eq{TS}, it follows that
\Eq{F2}{
&\int\limits_{T'_{x,y}}\int\limits_{T''_{x,y}}\Lambda(t',x,y)\Lambda(t'',x,y)\O(M_0(x,y),M(t'',x,
y))f(M(t',x,y))d\mu(t'')d\mu(t')\\
&=\big(\o_0(M_0(x,y))S_1''(x,y)-\o_1(M_0(x,y))S_0''(x,y)\big)
       \int\limits_{T'_{x,y}}\Lambda(t',x,y)f(M(t',x,y))d\mu(t')\\
&=\big((S_0'(x,y)+S_0''(x,y))S_1''(x,y)-(S_1'(x,y)+S_1''(x,y))S_0''(x,y)\big)
     \int\limits_{T'_{x,y}}\!\!\Lambda(t',x,y)f(M(t',x,y))d\mu(t')\\
&=\big(S_0'(x,y)S_1''(x,y)-S_1'(x,y)S_0''(x,y)\big)
     \int\limits_{T'_{x,y}}\!\!\Lambda(t',x,y)f(M(t',x,y))d\mu(t'),
}
and similarly,
\Eq{F3}{
\int\limits_{T'_{x,y}}\int\limits_{T''_{x,y}}\Lambda(t',x,y)\Lambda(t'',x,y)\O(M(t',x,y),M_0(x,
y))f(M(t'',x,y))d\mu(t'')d\mu(t')\\
=\big(S_1''(x,y)S_0'(x,y)-S_0''(x,y)S_1'(x,y)\big)
     \int\limits_{T''_{x,y}}\!\Lambda(t'',x,y)f(M(t'',x,y))d\mu(t'').
}
Substituting the above formulas \eq{F1}, \eq{F2}, and \eq{F3} into \eq{F0} and dividing the
inequality
so obtained by $(S_0'(x,y)S_1''(x,y)-S_1'(x,y)S_0''(x,y))$,
we get \eq{lH} with the error function $E:\Delta^\circ(I)\to\R$ defined by \eq{2E}, which completes
the proof.
\end{proof}

\Rem{1}{
A direct corollary of this theorem is the lower Hermite--Hadamard type inequality established by
\thm{E}.
Indeed, suppose that, with the notations introduced in  \eq{0kc}, the assumptions of \thm{E} hold.
Then, the $(\o_0,\o_2)$-convexity of $f$ implies that it is locally bounded and Borel measurable.
We show first that
the conditions of \thm{3} are also valid.
Let $\mu$ denote the Lebesgue measure on $[0,1]$ and define, for all $(x,y)\in \Delta^{\circ}(I)$,
$t\in[0,1]$,
\Eq{*}{
  M_0(x,y):=\xi(x,y),\quad M(t,x,y):=(1-t)x+ty, \quad\mbox{and}\quad
  \Lambda(t,x,y):=\frac{(y-x)\rho((1-t)x+ty)}{c(x,y)}.
}
Since $M(t,x,y)=M_0(x,y)$ can hold only for one value of $t$, hence \eq{2fi+} holds trivially.
We also have
\Eq{*}{
\int_0^1\Lambda(t,x,y)\o_1(M(t,x,y))dt&=\frac{y-x}{c(x,y)}\int_0^1\rho((1-t)x+ty)\o_1((1-t)x+ty)dt\\
                                      &= \frac{1}{c(x,y)}\int_x^y\rho\o_1
                                      =\o_0(\xi(x,y))\frac{\int_x^y\o_1\rho}{\int_x^y\o_0\rho}\\
                                      &=\o_0(\xi(x,y))\frac{\o_1}{\o_0}(\xi(x,y))
                                      =\o_1(\xi(x,y))=\o_1(M_0(x,y))
}
and, similarly,
\Eq{*}{
\int_0^1\Lambda(t,x,y)\o_0(M(t,x,y))dt=\frac{1}{c(x,y)}\int_x^y\rho\o_0
                                      =\o_0(\xi(x,y)),
}
which proves \eq{2=}. Thus all the assumptions (A1)--A(4) are verified. Therefore, if a function
$f:I\to\R$
is $(\o_0,\o_1)$-convex, i.e., satisfies \eq{J} with $\eps_{x,y}:=0$, then it fulfills \eq{lH} with
$E:=0$,
which, by the obvious identity $\frac{1}{c(x,y)}\int_x^yf\rho=\int_0^1\Lambda(t,x,y)f(M(t,x,y))dt$
is equivalent to
the left hand side inequality in \eq{EH}.
}

The following result could be deduced form \thm{3}, however a direct proof is more convenient here.
Given a set $D$, denote $\{(x,y)\mid x,y\in D,\,x\neq y\}$ by $D^{2*}$.

\Thm{4}{Let $D$ be a convex set of a linear space $X$. Let $\A$ be a sigma algebra containing the
Borel subsets of
$[0,1]$ and $\mu$ be a probability measure on the measure space $([0,1],\A)$ such that the support
of $\mu$ is
not a singleton. Denote
\Eq{*}{
   \mu_1:=\int\limits_{[0,1]}td\mu(t)\qquad\mbox{and}\qquad
S(\mu):=\mu\big([0,\mu_1]\big)\int\limits_{]\mu_1,1]}td\mu(t)-\mu\big(]\mu_1,1]\big)\int\limits_{[0,
\mu_1]}td\mu(t).
}
Assume that  $f:D\to \R$ is and hemi-$\mu$-integrable solution of the functional inequality
\Eq{JJ+}{
f((1-t)x+ty)\leq (1-t)f(x)+tf(y)+\eta_{x,y}(t)\qquad ((x,y)\in D^{2*},\,t\in [0,1])
}
where, for all $(x,y)\in D^{2*}$, $\eta_{x,y}:[0,1]\to \R$ is a function such that
\Eq{*}{
  I(x,y):=\int\limits_{]\mu_1,1]}\int\limits_{[0,\mu_1]}
          (t''-t')\eta_{(1-t')x+t'y,(1-t'')x+t''y}\Big(\frac{\mu_1-t'}{t''-t'}\Big)d\mu(t')d\mu(t'')
}
exists in $[-\infty,\infty]$ for all $(x,y)\in D^{2*}$.
Then, for all $(x,y)\in D^{2*}$, the function $f$ also satisfies
the lower Hermite--Hadamard type inequality
\Eq{lHc}{
  f((1-\mu_1) x+ \mu_1y)\leq \int\limits_{[0,1]}
f\big((1-t)x+ty)d\mu(t)+\frac{1}{S(\mu)}I(x,y)\qquad((x,y)\in D^{2*}).
}}

\Rem{T4}{In the above theorem, the hemi-$\mu$-integrability condition for $f$ can be relaxed if the
error function
$\eta_{x,y}$ enjoys boundedness and continuity properties. For instance, if $\eta_{x,y}$ is upper
bounded on $[x,y]$
for some $(x,y)\in D^{2*}$, then \eq{J} implies that $f((1-t)x+ty)$ is upper bounded for
$t\in[0,1]$. Similarly, if
$\limsup_{t\to 0+0}\eta_{x,y}(t)=0$ for some $(x,y)\in D^{2*}$, then \eq{J} implies that
$f((1-t)x+ty)$ is
an upper semicontinuous function of $t$ at zero from the right.}

\begin{proof} Let $(x,y)\in D^{2*}$.
Substituting in \eq{JJ+} $x$ by $(1-t')x+t'y$, $y$ by $(1-t'')x+t''y$ and $t$ by
$\frac{\mu_1-t'}{t''-t'}$,
where $0\leq t'\leq \mu_1$ and $\mu_1<t''\leq 1$, we get that
\Eq{4mu}{
f((1-\mu_1)x+\mu_1y)\leq
\frac{t''-\mu_1}{t''-t'}f((1-t')x&+t'y)+\frac{\mu_1-t'}{t''-t'}f((1-t'')x+t''y)\\
                    &+\eta_{(1-t')x+t'y,(1-t'')x+t''y}\Big(\frac{\mu_1-t'}{t''-t'}\Big).
}
Multiplying \eq{4mu} by $t''-t'$ and integrating on $[0,\mu_1]\times ]\mu_1,1]$ with respect to the
product measure
$\mu\times\mu$, we obtain
\Eq{F0+}{
\int\limits_{]\mu_1,1]}&\int\limits_{[0,\mu_1]}(t''-t')d\mu(t')d\mu(t'')f((1-\mu_1)x+\mu_1y)\\
                 &\leq
\int\limits_{]\mu_1,1]}(t''-\mu_1)d\mu(t'')\int\limits_{[0,\mu_1]}f((1-t')x+t'y)d\mu(t')\\          
&\quad+\int\limits_{[0,\mu_1]}(\mu_1-t')d\mu(t')\int\limits_{]\mu_1,1]}f((1-t'')x+t''y)d\mu(t'')\\
&\quad+\int\limits_{]\mu_1,1]}\int\limits_{[0,\mu_1]}(t''-t')\eta_{(1-t')x+t'y,(1-t'')x+t''y}
\Big(\frac{\mu_1-t'}{t''-t'}\Big)d\mu(t')d\mu(t'').
}
Applying Fubini's theorem, we get that
\Eq{F1+}{
\int\limits_{]\mu_1,1]}\int\limits_{[0,\mu_1]}(t''-t')d\mu(t')d\mu(t'')
=\mu\big([0,\mu_1]\big)\int\limits_{]\mu_1,1]}t''d\mu(t'')-\mu\big(]\mu_1,1]\big)\int\limits_{[0,
\mu_1]}t'd\mu(t')=S(\mu).
}
Using that the support of $\mu$ is not a singleton, we can see that the left hand side of \eq{F1+}
is positive
and hence so is $S(\mu)$.

Applying also Fubini's theorem, it follows that
\Eq{F2+}{
 \int\limits_{]\mu_1,1]}(&t''-\mu_1)d\mu(t'')
=\mu\big([0,1]\big)\int\limits_{]\mu_1,1]}t''d\mu(t'')-\mu\big(]\mu_1,1]\big)\int\limits_{[0,1]}
td\mu(t)\\
        &=\Big(\mu\big([0,\mu_1]\big)+\mu\big(]\mu_1,1]\big)\Big)\int\limits_{]\mu_1,1]}t''d\mu(t'')
-\mu\big(]\mu_1,1]\big)\Big(\int\limits_{[0,\mu_1]}t'd\mu(t')+\int\limits_{]\mu_1,1]}
t''d\mu(t'')\Big)=S(\mu)
}
and, similarly,
\Eq{F3+}{
\int\limits_{[0,\mu_1]}(\mu_1-t')d\mu(t')=S(\mu).
}
Substituting the above formulas \eq{F1+}, \eq{F2+}, and \eq{F3+} into \eq{F0+} and dividing the
inequality
so obtained by $S(\mu)$, we arrive at \eq{lHc}. This completes the proof.
\end{proof}

The following corollary is analogous to the result of \cite{HazPal09}.

\Cor{2HP1}{Assume that  $f:D\to \R$ a hemi-Lebesgue integrable solution of the functional
inequality \eq{JJ+}, where, for all $(x,y)\in D^{2*}$, $\eta_{x,y}:[0,1]\to \R$ is a function, such
that
\Eq{Ixy}{
  I(x,y):=\int\limits_{\frac12}^1\int\limits_{0}^{\frac12}(t''-t')\eta_{(1-t')x+t'y,(1-t'')x+t''y}
        \Big(\frac{\tfrac12-t'}{t''-t'}\Big)dt'dt''
}
exists in $[-\infty,\infty]$ for all $(x,y)\in D^{2*}$.
Then, for all $x,y\in D^{2*}$, the function $f$ also satisfies
\Eq{lHchh}{
  f\Big(\frac{x+y}{2}\Big)\leq \int\limits_0^1 f\big((1-t)x+ty)dt+8I(x,y).
}}

\begin{proof} We apply \thm{4}, when $\A$ is the family of Lebesgue measurable subsets of $[0,1]$,
$\mu$ is the Lebesgue measure. Then $\mu_1=\frac12$ and $S(\mu)=\frac18$ and the result directly
follows from \thm{4}.
\end{proof}

\Rem{2HP1}{In what follows, we deduce the conclusion of \thm{C} from the above corollary under
stronger
regularity assumption on $f$. Let $\alpha:(D-D)\to\R_+$ be a nonnegative
radially Lebesgue integrable function and assume that $f:D\to\R$ is hemi-upper bounded and
approximately
Jensen convex in the sense of \eq{0J}. Then, by \cite[Thm.\ 8]{MakPal11c}, $f$ fulfils the following
approximate convexity inequality:
\Eq{*}{
  f((1-t)x+ty)\leq (1-t)f(x)+tf(y)+\sum_{n=0}^\infty \frac1{2^n}\alpha(2d_\Z(2^nt)(x-y))
        \qquad ((x,y)\in D^{2},\,t\in [0,1]),
}
i.e., \eq{JJ+} holds with $\eta_{x,y}$ defined as
\Eq{*}{
  \eta_{x,y}(t):=\sum_{n=0}^\infty \frac1{2^n}\alpha(2d_\Z(2^nt)(x-y))\qquad((x,y)\in D^2,\,
t\in[0,1]).
}
Thus, by \cor{2HP1}, the inequality \eq{lHchh} holds with
\Eq{*}{
 I(x,y)&=\int\limits_{\frac12}^1\int\limits_{0}^{\frac12}(t''-t')\eta_{(1-t')x+t'y,(1-t'')x+t''y}
        \Big(\frac{\tfrac12-t'}{t''-t'}\Big)dt'dt''\\
  &=\sum_{n=0}^\infty \frac1{2^n}\int\limits_{\frac12}^1\int\limits_{0}^{\frac12}(t''-t')
        \alpha\Big(2d_\Z\Big(2^n\frac{\tfrac12-t'}{t''-t'}\Big)(t''-t')(x-y)\Big)dt'dt''\\
  &=\sum_{n=0}^\infty \frac1{2^n}\int\limits_0^{\frac12}\int\limits_{0}^{\frac12}(t+s)
        \alpha\Big(2d_\Z\Big(\frac{2^nt}{t+s}\Big)(t+s)(x-y)\Big)dtds\\
  &=\sum_{n=0}^\infty \frac2{2^n}\int\limits_0^{\frac12}\int\limits_{0}^{t}(t+s)
        \alpha\Big(2d_\Z\Big(\frac{2^nt}{t+s}\Big)(t+s)(x-y)\Big)dsdt.
}
The last equality above is the consequence of the symmetry of the integrand with respect to the
variables $s,t$. For $n=0$,
\Eq{*}{
\frac2{2^n}\int\limits_0^{\frac12}\int\limits_{0}^{t}(t+s)
        &\alpha\Big(2d_\Z\Big(\frac{2^nt}{t+s}\Big)(t+s)(x-y)\Big)dsdt
  =2\int\limits_0^{\frac12}\int\limits_{0}^{t}(t+s)\alpha(2s(x-y))dsdt\\
 &=2\int\limits_0^{\frac12}\int\limits_{s}^{\frac12}(t+s)\alpha(2s(x-y))dtds
  =\int\limits_0^{\frac12}(1-2s)\Big(\frac32s+\frac14\Big)\alpha(2s(x-y))ds\\
 &=\frac18\int\limits_0^1(1-\sigma)(3\sigma+1)\alpha(\sigma(x-y))d\sigma.
}
To compute the the double integral for $n\geq1$, we will split its domain according to the position
of $\frac{2^nt}{t+s}$ related to integer numbers. For all $n\in \N$ and $0<s<t\leq \frac12$, there
exists a unique $m\in\{2^{n-1},\dots,2^n-1\}$ (namely $m:=\left[\frac{2^nt}{t+s}\right]$) such that
\Eq{*}{
\mbox{either}\qquad m\leq \frac{2^nt}{t+s}< m+\frac12\qquad \mbox{or}\qquad
m+\frac12\leq\frac{2^nt}{t+s}<m+1.
}
This, for all $m\in\{2^{n-1},\dots,2^n-1\}$, in terms of $t$ yields the following inequalities for
$s$:
\Eq{*}{
\frac{2^n-m-\frac12}{m+\frac12}t < s \leq \frac{2^n-m}{m}t \qquad \mbox{and}\qquad
\frac{2^n-m-1}{m+1}t < s \leq \frac{2^n-m-\frac12}{m+\frac12}t,
}
respectively. On these intervals, we have that
\Eq{*}{
  d_\Z\Big(\frac{2^nt}{t+s}\Big)(t+s)=
  \begin{cases}
   (\frac{2^nt}{t+s}-m)(t+s)=(2^n-m)t-ms,
            \quad \mbox{if}\quad\frac{2^n-m-\frac12}{m+\frac12}t < s \leq \frac{2^n-m}{m}t, \\[2mm]
   (m+1-\frac{2^nt}{t+s})(t+s)\!=\!(m+1-2^n)t+(m+1)s,
            \quad \!\!\mbox{if}\!\!\quad\frac{2^n-m-1}{m+1}t < s \leq
\frac{2^n-m-\frac12}{m+\frac12}t.
  \end{cases}
}
Thus, we get that
\Eq{*}{
&\int\limits_0^{\frac12}\int\limits_{0}^{t}(t+s)
        \alpha\Big(2d_\Z\Big(2^n\frac{t}{t+s}\Big)(t+s)(x-y)\Big)dsdt\\
        &=\int\limits_0^{\frac12}
     \sum_{m=2^{n-1}}^{2^n-1}\bigg(\int\limits_{\frac{2^n-m-\frac12}{m+\frac12}t}^{\frac{2^n-m}{m}t}
             (t+s)\alpha\big(2((2^n-m)t-ms)(x-y)\big)ds\\
      &\qquad\qquad\qquad+\int\limits_{\frac{2^n-m-1}{m-1}t}^{\frac{2^n-m-\frac12}{m+\frac12}t}(t+s)
        \alpha\big(2((m+1-2^n)t+(m+1)s)(x-y)\big)ds\bigg)dt\\
        &=\sum_{m=2^{n-1}}^{2^n-1}\int\limits_0^{\frac12}
\bigg(\int\limits_{0}^{\frac{2^{n+1}t}{2m+1}}
            \alpha(\sigma(x-y))\Big(\frac{\sigma+2^{n+1}t}{(2m+2)^2}
        +\frac{2^{n+1}t-\sigma}{(2m)^2}\Big)d\sigma\bigg)dt\\
        &=\sum_{m=2^{n-1}}^{2^n-1}\int\limits_0^{\frac{2^{n}}{2m+1}}
            \bigg(\int\limits_{\frac{(2m+1)\sigma}{2^{n+1}}}^{\frac12}
            \alpha(\sigma(x-y))\Big(\frac{\sigma+2^{n+1}t}{(2m+2)^2}
        +\frac{2^{n+1}t-\sigma}{(2m)^2}\Big)dt\bigg)d\sigma\\
        &=\frac1{16}\sum_{m=2^{n-1}}^{2^n-1}\int\limits_0^{\frac{2^{n}}{2m+1}}
            \alpha(\sigma(x-y))\Big(1-\frac{2m+1}{2^n}\sigma\Big)
            \Big(\frac{\sigma(2m+3)+2^{n}}{(m+1)^2}
        +\frac{\sigma(2m-1)+2^n}{m^2}\Big)d\sigma\\
        &=\frac1{16}\sum_{m=2^{n-1}}^{2^n-1}\int\limits_0^1
            \alpha(\sigma(x-y))\Big(1-\frac{2m+1}{2^n}\sigma\Big)^+
            \Big(\frac{\sigma(2m+3)+2^{n}}{(m+1)^2}
        +\frac{\sigma(2m-1)+2^n}{m^2}\Big)d\sigma.
}
(Here $x^+$ stands for the positive part of $x$.) Summarizing our computations, for $8I(x,y)$, we
get
\Eq{*}{
  8I(x,y)=\int\limits_0^1&\alpha(\sigma(x-y))\Phi(\sigma)d\sigma,
}
where
\Eq{*}{
  \Phi(\sigma):
  &=(1-\sigma)(3\sigma+1)+\sum_{n=1}^\infty\sum_{m=2^{n-1}}^{2^n-1}
            \Big(1-\frac{2m+1}{2^n}\sigma\Big)^+\Big(\frac{\sigma(2m+3)+2^{n}}{2^n(m+1)^2}
        +\frac{\sigma(2m-1)+2^n}{2^nm^2}\Big)\\
  &=(1-\sigma)(3\sigma+1)+\sum_{m=1}^\infty
\Big(1-\frac{2m+1}{2^{[\log_2m]+1}}\sigma\Big)^
+\Big(\frac{\sigma(2m+3)+2^{[\log_2m]+1}}{2^{[\log_2m]+1}(m+1)^2}
        +\frac{\sigma(2m-1)+2^{[\log_2m]+1}}{2^{[\log_2m]+1}m^2}\Big).
}
One can easily see that $\Phi$ is a continuous function over $[0,1]$ with $\Phi(t)\geq1$ for $0\leq
t\leq\frac23$
and $\Phi(1)=0$. Hence the error term $8I(x,y)$ obtained in \eq{lHchh} is not comparable with that
in \eq{0lH}.}

In what follows, we examine the case, when $X$ is a normed space and
$\eta_{x,y}(t)$ is a linear combination of the products of the powers of $t$, $1-t$, and
of $\|x-y\|$, i.e., for all $(x,y)\in D^{2*}$ $\eta_{x,y}$ is of the form
\Eq{4DHP}{
  \eta_{x,y}(t):=\int\limits_{[0,\infty[^2} t^p(1-t)^q\|x-y\|^{p+q-1}d\nu(p,q)\qquad ((x,y)\in
D^{2*}),
}
where $\nu$ is a $\sigma$-finite Borel measure on $[0,\infty[^2$.
An important particular case is when $\nu$ is of the form $\sum_{i=1}^k c_i\delta_{(p_i,q_i)}$,
where $c_1,\dots,c_k>0$, $(p_1,q_1),\dots,(p_k,q_k)\in[0,\infty[^2$.

\Thm{3A2}{Let $\A$ be a sigma algebra containing the Borel subsets of $[0,1]$ and $\mu$ be a
probability
measure on the measure space $([0,1],\A)$ such that the support of $\mu$ is not a singleton.
Let $\nu$ be a $\sigma$-finite Borel measure on $[0,\infty[^2$ such that, for all
$s\in\{\|x-y\|\mid (x,y)\in D^{2*}\}$,
\Eq{*}{
  J(s):=\int\limits_{[0,\infty[^2}\Big(\int\limits_{[0,\mu_1]}(\mu_1-t')^pd\mu(t')
     \int\limits_{]\mu_1,1]}(t''-\mu_1)^qd\mu(t'')\Big)s^{p+q-1}d\nu(p,q)
}
exists in $[-\infty,\infty]$.
Assume that  $f:D\to \R$ is a hemi-$\mu$-integrable solution of the functional inequality
\Eq{3Jnu+}{
  f((1-t)x+ty)\leq (1-t)f(x)+tf(y)
        +\int\limits_{[0,\infty[^2} t^p(1-t)^q\|x-y\|^{p+q-1}d\nu(p,q)
}
for all $(x,y)\in D^{2*}$ and $t\in[0,1]$.
Then $f$ also fulfils the Hermite--Hadamard type inequality
\Eq{lHc+}{
  f((1-\mu_1) x+ \mu_1y)\leq \int\limits_{[0,1]} f\big((1-t)x+ty)d\mu(t)+\frac{1}{S(\mu)}J(\|x-y\|)
  \qquad((x,y)\in D^{2*}).
}}

\Rem{T3A2}{In the above theorem, the hemi-$\mu$-integrability condition for $f$ can be relaxed if
the measure $\nu$ is finite with compact support contained in $]0,\infty[^2$. Then the function
$\eta_{x,y}$ defined by \eq{4DHP} is continuous on $[x,y]$ and $\eta_{x,y}(0)=\eta_{x,y}(1)=0$,
hence \eq{3Jnu+} implies that $t\mapsto f((1-t)x+ty)$ is upper bounded on $[0,1]$ and upper
semicontinuous at the endpoint of $[0,1]$. Thus $f$ is hemi-upper bounded and
upper hemicontinuous on $D$, which yields its hemi-$\mu$-integrability.}

\begin{proof} We want to apply \thm{4}. Let $(x,y)\in D^{2*}$ be arbitrary. Let
$\eta_{x,y}:[0,1]\to\R$
defined by \eq{4DHP}. Then, \eq{3Jnu+} is equivalent to \eq{JJ+}. To deduce \eq{lHc+}, by \thm{4},
we obtain that
\Eq{*}{
I(x,y)
&=\int\limits_{]\mu_1,1]}\int\limits_{[0,\mu_1]}(t''-t')
    \int\limits_{[0,\infty[^2}\big(\tfrac{\mu_1-t'}{t''-t'}\big)^p
   \big(\tfrac{t''-\mu_1}{t''-t'}\big)^q\|(t''-t')(x-y)\|^{p+q-1}d\nu(p,q)d\mu(t')d\mu(t'')\\
&=\int\limits_{[0,\infty[^2}\Big(\int\limits_{[0,\mu_1]}(\mu_1-t')^pd\mu(t')
     \int\limits_{]\mu_1,1]}(t''-\mu_1)^qd\mu(t'')\Big)\|x-y\|^{p+q-1}d\nu(p,q)=J(\|x-y\|),
}
which proves the statement.
\end{proof}

\Cor{4C2}{Let $\nu$ be a $\sigma$-finite Borel measure on $[0,\infty[^2$, such that
for all $s \in \{\|x-y\|:(x,y)\in D^{2*}\}$
\Eq{*}{
\int\limits_{[0,\infty[^2}\frac{s^{p+q-1}}{2^{p+q-1}(p+1)(q+1)}d\nu(p,q)
}
exists in $[-\infty,\infty]$.
Assume that $f:D\to \R$ is a hemi-Lebesgue integrable solution of the functional inequality
\Eq{C2a}{
f\big((1-t)x+ty\big)\leq (1-t)f(x)+tf(y)
+\int\limits_{[0,\infty[^2}t^p(1-t)^q\|x-y\|^{p+q-1}d\nu(p,q),
}
where $(x,y)\in D^{2*}$ and $t\in[0,1]$.
Then, $f$ also satisfies the Hermite--Hadamard type inequality
\Eq{C2b}{
  f\Big(\frac{x+y}{2}\Big)\leq \int\limits_{[0,1]} f\big((1-t)x+ty)dt
    +\int\limits_{[0,\infty[^2}\frac{\|x-y\|^{p+q-1}}{2^{p+q-1}(p+1)(q+1)}d\nu(p,q) \qquad((x,y)\in
D^{2*}).
}}

\begin{proof}Observe that \eq{C2a} is equivalent to \eq{JJ+}, where for all $(x,y)\in D^{2*}$, 
$\eta_{x,y}:[0,1]\to\R$ is defined by \eq{4DHP}. We have $S(\mu)=\frac18$ and using \thm{3A2}, we
obtain that
\Eq{*}{
J&(s)=\int\limits_{[0,\infty[^2}
       \int\limits_0^{\frac12}(\tfrac12-t)^pdt\int\limits_{\frac12}^1(t-\tfrac12)^qdt
s^{p+q-1}d\nu(p,q)\!
        =\!\!\int\limits_{[0,\infty[^2}\frac{s^{p+q-1}}{2^{p+q+2}(p+1)(q+1)}d\nu(p,q),
}
which yields \eq{C2b}.
\end{proof}

\section{From approximate $(\o_0,\o_1)$-convexity to approximate upper Hermite--Hadamard inequality}
\setcounter{equation}{0}

In the first part of this section we will investigate the implication between the
$(\o_0,\o_1)$-convexity
type inequality and upper Hermite--Hadamard inequality. Consider the following assumptions.
\begin{enumerate}[({B}1)]
\item $(T,\A,\mu)$ is a measure space.
\item $\Lambda:T\times \Delta^\circ(I)\to\R_+$ is integrable (with respect to $\mu$) in its first
variable.
\item $M:T\times \Delta^\circ(I)\to\R$ is measurable in its first variable
and for all $t\in T$, the map $(x,y)\mapsto M(t,x,y)$ is a two-variable mean on $I$.
$M_0:\Delta^\circ(I)\to I$ is a strict mean. 
\item There exist an $(\o_0,\o_1)$-Chebyshev system on
$I$ such that $\o_0$ is positive and $i\in\{0,1\}$ \eq{2=} holds.
\end{enumerate}

\Thm{5}{Assume that (B1)--(B4) hold. Let $f:I\to\R$ be a locally bounded Borel measurable solution
of the approximate $(\o_0,\o_1)$-convexity inequality \eq{J}, where for all
$(x,y)\in\Delta^{\circ}(I)$, $\eta_{x,y}:[x,y]\to\R$ is a bounded and Borel measurable function.
Then $f$ also satisfies the following approximate upper Hermite--Hadamard type inequality
\Eq{uH}{
\int\limits_{T}\Lambda(t,x,y)f(M(t,x,y))d\mu(t)
\leq \frac{\O(M_0(x,y),y)}{\O(x,y)}f(x)+\frac{\O(x,M_0(x,y))}{\O(x,y)}+E(x,y),
}
with $E:D^{2*}\to\R$ defined by
\Eq{uE}{
E(x,y)=\int\limits_T\Lambda(t,x,y)\eps_{x,y}(M(t,x,y))d\mu(t).
}
}

\begin{proof}Let $(x,y)\in \Delta^{\circ}(I)$ be arbitrary. Substituting in \eq{J} $u$ by
$M(t,x,y)$, we get that 
\Eq{*}{
f(M(t,x,y))\leq
\frac{\O(M(t,x,y),y)}{\O(x,y)}f(x)+\frac{\O(x,M(t,x,y))}{\O(x,y)}f(y)+\eps_{x,y}(M(t,x,y))
\qquad(t\in T).
}
Multiplying this inequality by $\Lambda(t,x,y)$ and integrating with respect to $\mu$ on $T$, we get
that 
\Eq{G0}{
\int\limits_T\Lambda(t,x,y)&f(M(t,x,y))d\mu(t)\\
&\leq\frac{\int\limits_T\Lambda(t,x,y)\O(M(t,x,y),y)d\mu(t)}{\O(x,y)}f(x)
  +\frac{\int\limits_T\Lambda(t,x,y)\O(x,M(t,x,y))d\mu(t)}{\O(x,y)}f(y)\\
  &+\int\limits_T\Lambda(t,x,y)\eps_{x,y}(M(t,x,y))d\mu(t).
}
Applying \eq{2=}, it follows that
\Eq{G1}{
\int\limits_T\Lambda(t,x,y)\O(M(t,x,y),y)d\mu(t)=\O(M_0(x,y),y)
}
and
\Eq{G2}{
\int\limits_T\Lambda(t,x,y)\O(x,M(t,x,y))d\mu(t)=\O(x,M_0(x,y)).
}
Substituting \eq{G1} and \eq{G2} to \eq{G0} we have \eq{uH}.
\end{proof}

\Rem{3}{An immediate corollary of this theorem is the second inequality of \thm{E}. Assume that the
assumptions of \thm{E} hold. It is easy to see that the conditions of \thm{5} are also valid. For
all $(x,y)\in \Delta^{\circ}(I)$ let $\eta_{x,y}$,$\mu$, $M(t,x,y)$, $M_0(x,y)$ and $\Lambda(t,x,y)$
be defined as in \rem{1}. Then \eq{2=} holds. Therefore, by \eq{0kc} and using also \rem{1},
\Eq{*}{
c_1(x,y)=\frac{1}{\O(x,y)}\left| \begin{array}{ccc}
\int_x^y\o_0\rho &\o_0(y) \\[2mm]
\int_x^y\o_1\rho & \o_1(y)\end{array} \right|
=\frac{1}{\O(x,y)}\left| \begin{array}{ccc}
\o_0(M_0(x,y)) &\o_0(y) \\[2mm]
\o_1(M_0(x,y)) & \o_1(y)\end{array} \right|=\frac{\O(M_0(x,y),y)}{\O(x,y)}.
}
Similarly, it can be seen, that
$c_2(x,y)=\frac{\O(x,M_0(x,y))}{\O(x,y)}$. Thus, by \thm{5}, we get the second inequality in
\thm{E}.
}

\Thm{6}{Let $D$ be a convex set of a linear space $X$. Let $\A$ be a sigma algebra containing the
Borel subsets of $[0,1]$ and $\mu$ be a probability measure on the measure space $([0,1],\A)$.
Denote $\mu_1:=\int_{[0,1]}td\mu(t)$. Assume that  $f:D\to \R$ is a hemi-$\mu$-integrable solution
of the approximate convexity inequality \eq{JJ+}, where, for all $(x,y)\in D^{2*}$,
$\eta_{x,y}:[0,1]\to \R$ is a bounded and Borel measurable function. Then, for all $(x,y)\in
D^{2*}$, the function $f$ also satisfies the approximate upper Hermite--Hadamard inequality
\Eq{uHc}{
\int\limits_{[0,1]} f\big((1-t)x+ty)d\mu(t)\leq
(1-\mu_1)f(x)+\mu_1f(y)+\int\limits_{[0,1]}\eta_{x,y}(t)d\mu(t).
}}

\Rem{T6}{In the above theorem, the regularity condition for $f$ can be relaxed if the error
function $\eta_{x,y}$ enjoys boundedness and continuity properties. For instance, if $\eta_{x,y}$ is
upper bounded on $[x,y]$ for some $(x,y)\in D^{2*}$, then \eq{J} implies that $f((1-t)x+ty)$ is
upper bounded for $t\in[0,1]$. Similarly, if $\limsup_{t\to 0+0}\eta_{x,y}(t)=0$ for some $(x,y)\in
D^{2*}$, then \eq{J} implies that $f((1-t)x+ty)$ is an upper semicontinuous function of $t$ at zero
from the right.}

\begin{proof}Let $(x,y)\in D^{2*}$ be fixed. Integrating \eq{JJ+} with respect to the variable $t$
and the measure $\mu$ on $[0,1]$ we get \eq{uHc}.
\end{proof}

\Rem{4}{
Assume that the conditions of \thm{A} hold. To prove a similar result as in \thm{A}, we have to
assume that also $\alpha:(D-D)\to\R$ is radially bounded and radially continuous at $0$. By
\cite{TabTab09b} and \cite{MakPal11c}, we can get that $f$ is approximately convex in the following
sense
\Eq{*}{
f((1-t)x+ty)\leq (1-t)f(x)+tf(y)+\sum_{n=0}^{\infty}\frac{1}{2^n}\alpha(2d_\Z(2^nt)(x-y))
 \qquad (x,y\in D,\, t\in[0,1]).
}
Let $\eta_{x,y}(t):=\sum_{n=0}^{\infty}\frac{1}{2^n}\alpha(2d_\Z(2^nt)(x-y))$ for $t\in[0,1]$ and
$x,y\in D$. Let $\A$ be the class of Lebesgue measurable subsets of $[0,1]$ and let the measure
$\mu$ be defined by $d\mu(t)=\rho(t)dt$. Then $\mu_1=\int_0^1t\rho(t)dt=\lambda$. Thus applying
\thm{6} and the Fubini's theorem, we get \eq{0H}, which completes the proof of \thm{A}.
}

In what follows, we examine the case, when $X$ is a normed space and
$\eta_{x,y}(t)$ is a linear combination of the products of the powers of $t$, $1-t$, and
of $\|x-y\|$, i.e., for all $(x,y)\in D^{2*}$ $\eta_{x,y}$ is of the form
\Eq{5DHP}{
  \eta_{x,y}(t):=\int\limits_{[0,\infty[^3} t^p(1-t)^q\|x-y\|^{r}d\nu(p,q,r)\qquad ((x,y)\in
D^{2*}),
}
where $\nu$ is a $\sigma$-finite Borel measure on $[0,\infty[^3$.
An important particular case is when $\nu$ is of the form $\sum_{i=1}^k c_i\delta_{(p_i,q_i,r_i)}$,
where $c_1,\dots,c_k>0$ and $(p_1,q_1,r_1),\dots,(p_k,q_k,r_k)\in[0,\infty[^3$.

\Cor{6a}{Let $\mu$ be a Borel probability measure on $[0,1]$, denote $\mu_1:=\int_{[0,1]}td\mu(t)$.
Let $\nu$ be a $\sigma$-finite Borel measure on $[0,\infty[^3$ such that, for all
$s\in\{\|x-y\|\mid (x,y)\in D^{2*}\}$,
\Eq{*}{
 \int\limits_{[0,\infty[^3} \int\limits_{[0,1]}t^p(1-t)^qs^{r}d\mu(t)d\nu(p,q,r)
}
exists in $[-\infty,\infty]$.
Assume that  $f:D\to \R$ is  and hemi-$\mu$-integrable solution of the functional inequality
\Eq{6Jnu+}{
  f((1-t)x+ty)\leq (1-t)f(x)+tf(y)
        +\int\limits_{[0,\infty[^3}  t^p(1-t)^q\|x-y\|^{r}d\nu(p,q,r)
}
for all $(x,y)\in D^{2*}$ and $t\in[0,1]$.
Then $f$ also fulfills the following approximate upper Hermite--Hadamard inequality,
\Eq{uH2}{
\int\limits_{[0,1]} f\big((1-t)x+ty)d\mu(t)\leq (1-\mu_1)f(x)+\mu_1f(y)
+\int\limits_{[0,\infty[^3}\int\limits_{[0,1]}t^p(1-t)^q\|x-y\|^{r}d\mu(t)d\nu(p,q,r).
}}
\begin{proof}
We apply \thm{6}. For all $x,y\in D^{2*}$, let $\eta_{x,y}:[0,1]\to\R$ defined by
\eq{5DHP}. Then it is easy to see that \eq{6Jnu+} is equivalent to \eq{JJ+}.
Hence, by \thm{6} we get \eq{uH2}.
\end{proof}

Denote by $B$ the so-called beta-function, defined by
\Eq{*}{
B(p_1,p_2)=\int\limits_0^1t^{p_1-1}(1-t)^{p_2-1}dt \qquad (p_1,p_2>0).
}

\Cor{6b}{Let $\nu$ be a $\sigma$-finite Borel measure on $[0,\infty[^3$ such that, for all
$s\in\{\|x-y\|\mid (x,y)\in D^{2*}\}$,
\Eq{*}{
 \int\limits_{[0,\infty[^3}B(p+1,q+1)s^{r}d\nu(p,q,r)
}
exists in $[-\infty,\infty]$.
Assume that  $f:D\to \R$ is a hemi-Lebesgue integrable solution of the functional inequality
\Eq{*}{
  f((1-t)x+ty)\leq (1-t)f(x)+tf(y)
        +\int\limits_{[0,\infty[^3}  t^p(1-t)^q\|x-y\|^{r}d\nu(p,q,r)
}
for all $(x,y)\in D^{2*}$ and $t\in[0,1]$.
Then $f$ also fulfills the approximate upper Hermite--Hadamard inequality
\Eq{uH3}{
f\Big(\frac{x+y}{2}\Big)\leq
\frac{f(x)+f(y)}{2}+\int\limits_{[0,\infty[^3}B(p+1,q+1)\|x-y\|^{r}d\nu(p,q,r)
    \qquad ((x,y)\in D^{2*}).
}}

\begin{proof}
We apply \cor{6a} when $\mu$ is the Lebesgue measure. Then, for all $(x,y)\in D^{2*}$,
\Eq{*}{
E(x,y)=\int\limits_{[0,\infty[^3}\int\limits_{[0,1]}t^p(1-t)^q\|x-y\|^{r}d\mu(t)d\nu(p,q,r)
      =\int\limits_{[0,\infty[^3}B(p+1,q+1)\|x-y\|^{r}d\nu(p,q,r).
}
Thus, the result directly follows from \cor{6a}.
\end{proof}


\end{document}